\documentclass[11pt]{article}

\usepackage{publication_style}
\hypersetup{pdftitle={Counting Classical Fullerenes VII: Arithmetic Unranking of Fullerene Isomers}}
\title{Counting Classical Fullerenes VII:\\Arithmetic Unranking of Fullerene Isomers}
\author{Kerem Kat\\[2pt]\small Independent Researcher\\\small\texttt{keremkat@gmail.com}}
\date{}
\begin{document}
\maketitle

\begin{abstract}
Given a vertex number $2m$ and a rank, we construct the corresponding
ordinary fullerene as a planar adjacency list. Choosing the rank uniformly
therefore samples ordinary isomorphism classes exactly from the vertex
number alone. We select an integral period by counting lattice points
on a quadratic level in polytopes with fixed facet directions.
For a fixed positive semidefinite integral form on a box of radius $H$,
with a linear term of coefficients $O(H)$, the algorithm computes all
level counts simultaneously in $O(H^2\operatorname{polylog}(2m))$
bit operations when $H=O(\sqrt m)$ and the facet offsets have
$O(\log(2m))$ bits. A signed decomposition into orthogonal simplices
gives this bound. For fullerenes we use deterministic
$O(m^2\operatorname{polylog}(2m))$ bit work to construct an index
of $O(m\log(2m))$ bits. With this index a rank costs
$O(m\operatorname{polylog}(2m))$ bit work, uniformly in the rank.
Two families of recuts let energy descent escape self-contacts in the
slit polygon. The decoder then writes the graph from its selected
period in $O(m)$ arithmetic word operations.
\end{abstract}

\section{Introduction}
A classical fullerene is a cubic three-connected planar graph with
only pentagonal and hexagonal faces. Write $\iso(2m)$ for the number
of ordinary isomorphism classes on $2m$ vertices, identifying mirror
images. A dense unranking assigns the integers $0,\ldots,\iso(2m)-1$
bijectively to the classes and constructs the graph at any given rank.

Given only a feasible vertex count, how can we draw each ordinary
isomer with the same probability? Bille, Buchstaber and Spodarev
ask for an algorithm that neither enumerates the isomers nor requires
their cardinality as input \cite[Sec.~3, Open Problem~1]{bille2025open}.
We compute the size of the rank interval and construct the graph at
a uniformly chosen rank. Each class has one rank, including those
with nontrivial automorphisms (\cref{cor:uniform-sampling}).

The counting step applies beyond fullerenes. Fix a dimension $d$, a
positive semidefinite integral quadratic form $S$ and a finite set of
facet directions. In a polytope inside $[-CH,CH]^d$, consider the lattice
points satisfying $S(z)+L(z)=n$, where the coefficients of the linear
form $L$ are $O(H)$. For $H=O(\sqrt m)$ and logarithmic-size input
integers, \cref{lem:fixed-normal-orthoscheme} counts these points in
$O(H^2\operatorname{polylog}(2m))$ bit operations. The exponent of $H$
is independent of $d$. The same computation gives the entire table of
level counts, which can have $H^2$ nonzero entries. A cut on the next coordinate preserves the
hypotheses, so the same bound gives random access in lexicographic order.

\begin{theorem}[Arithmetic unranking]\label{thm:main}
For every positive integer $m$ there is an explicitly defined
lexicographic order of the ordinary fullerenes on $2m$ vertices.
In deterministic $O(m^2\operatorname{polylog}(2m))$ bit operations
one can compute their count and build an index occupying
$O(m\log(2m))$ bits. Given this index, any valid zero-based rank
returns a graph in $O(m\operatorname{polylog}(2m))$ bit operations.
Both procedures use $O(m\operatorname{polylog}(2m))$ bits of working
storage. Once its period is selected, the graph takes $O(m)$
arithmetic operations on $O(\log(2m))$-bit words, including output.
\end{theorem}
Choose one period representative per isomer to obtain a unique code
(\cref{sec:owners}). Codes sharing a prefix occupy consecutive ranks.
Counting completions for each next-coordinate threshold locates the
block containing the requested rank. Subtract the size of the preceding
blocks and repeat until the period is determined
(\cref{thm:orthoscheme-family-count}). Recut its slit polygon until
the boundary is simple, glue the triangular mesh, and read off its
dual graph (\cref{thm:star-complete-recut-normalization,sec:decoder}).

A planar adjacency list records the three neighbours of each vertex
in cyclic order, with arbitrary vertex labels.
We count ordinary bit operations, with logarithmic factors allowed
for exact integer arithmetic. In the separate word bound we allow
integer multiplication and division on logarithmic-size words.
Constants depend on a single fixed lattice of Eisenstein rank ten.
The bounds concern a mathematical algorithm and give no crossover
estimate against graph generation. The count also detects an empty
family or a rank outside the admissible interval.
For an explicit list on vertices labelled $1,\ldots,2m$, cubicity
forces each label to occur three times among the neighbours. Writing
those entries takes $\Omega(m)$ words, or $\Omega(m\log(2m))$ bits
when the labels are written in binary. This gives an output lower
bound against which to compare the indexed query time.

\section{Prior work on generation and reconstruction}
Recursive generation gives a complete census by local graph
expansions \cite[Secs.~2--3]{brinkmann2012}. Graver reconstructs
a fullerene from its labelled twelve-vertex signature
\cite[Theorem~1 and Sec.~6]{graver2004encoding}.
Wirz, Schwerdtfeger and Avery give a canonical general face-spiral
name and a windup algorithm that reconstructs its polyhedron
\cite[Sec.~II.E, Algorithm~2, and Sec.~III]{wirz2018faceSpirals}.
Both constructions take an encoding of the individual graph as input.
Thurston's lattice correspondence identifies convex triangulated cone
spheres through their periods \cite[Theorem~0.1 and Sec.~7]{thurston1998}.
Here a size and a rank determine the period
by counts of completions, so selection does not enumerate earlier graphs.

\paragraph{Period conventions.}
Let $\mathcal E=\ZZ[\zeta]$, $\zeta=e^{\pi i/3}$ and $\theta=1+\zeta$.
In a basis $\alpha_0,\ldots,\alpha_9$ the Hermitian lattice
$\Lambda=\mathcal E^{10}$ has Gram matrix with diagonal $-3$, upper
adjacent entry $\theta$, its conjugate below, and all other entries zero.
The form $h$ has signature $(1,9)$ and is linear in the second variable.
A root has norm $-3$. A vector is open, or mirror-free, if its
perpendicular contains no root. By \cref{cnt-prop:fullerene-period-correspondence} of~\cite{fullereneCountManuscript}
and its orientation-reversal description, ordinary isomers correspond
to orbits of open vectors of norm $3m$ under the group generated by
integral unitary isometries and an integral antiunitary lift of
orientation reversal. Taking a scalar unit
does not change the orbit. Imprimitive vectors are included. This is
the sole geometric counting input from Paper~I.

\section{One integral tuple per ordinary orbit}\label{sec:owners}
The first problem is to retain exactly one period in each orbit
using inequalities suitable for counting. A uniform bound for the
smallest cusp height reduces the search to finitely many coefficients.

\subsection{Integral cusp coordinates}
Write $\mathcal E=\mathbb Z[\zeta]$, $\zeta=e^{\pi i/3}$ and
$\theta=1+\zeta$.  The Hermitian form is linear in its second variable.
At a primitive cusp use
\[
 \Lambda=(\mathcal E e+\mathcal E f)\perp F(-1)\perp F(-1),
 \qquad h(e,f)=\bar\theta,\qquad F=E_8^{\mathcal E}.
\]
All primitive cusps have this type: the twelve equal curvatures have the
single cusp partition $6+6$ \cite[Section~5]{thurston1998}, and the
primitive-null-vector correspondence is
\cite[Lemma~2]{allcock2008eisenstein}.  This is also the splitting in
\cref{ipr-lem:ipr-cusp-fiber} of~\cite{fullereneIprManuscript}.  On the positive horizontal lattice
$K=F\perp F$ put $Q(x)=h_+(x,x)/3$.  Its underlying real lattice is
$E_8\perp E_8$ and $Q$ is half its squared Euclidean norm.

Here is one convention fixing all coordinate choices for the order.
In the star chart of \cref{prop:star-marked-period-chart} choose
$e$ from the six successive edges $1,\zeta,\ldots,\zeta^5$ followed
by five zero edges and anchored at the origin. Their tips give a
primitive null vector by \eqref{eq:star-area-shell}. Among integral
bases $(e,f,k_1,\ldots,k_8)$ require $f$ null,
$h(e,f)=\bar\theta$ and the $k_i$ orthogonal to $e,f$. On the
$k_i$ use two mutually orthogonal four-node chains for $h_+$, each
with diagonal $3$, upper adjacent entry $-\bar\theta$ and all other
upper entries zero. Choose the first such basis in increasing order of the largest absolute integer
coefficient and then in lexicographic order of all coefficients
in the $\alpha$ basis, writing the $1$ coefficient before the
$\zeta$ coefficient. The splitting ensures the search terminates.
All coordinate comparisons below use this ordered basis of $K$.

For every $t\in K$ the following integral isometry fixes $e$:
\begin{equation}\label{eq:cusp-inventory-translation}
 T_t(f)=f+t+Q(t)e,\qquad
 T_t(x)=x+\frac{h_+(t,x)}{\theta}e\quad(x\in K).
\end{equation}
Here $Q(t)\in\mathbb Z$ and $h_+(K,K)\subseteq\theta\mathcal E$.
For example,
$h(T_t(x),T_t(f))=-h_+(x,t)+h_+(x,t)=0$, and
$h(T_t(f),T_t(f))=3Q(t)-h_+(t,t)=0$.
The other pairings follow directly; the matrix is integral and
triangular with determinant one.  The central isometries are
$f\mapsto f+k\zeta e$ for $k\in\mathbb Z$.
Consequently, for $v=u e+\beta f+z$, horizontal translations change
$z$ by $\beta t$ and central translations change $u$ by $k\beta\zeta$.
These are elements of the full unitary group, without a congruence
restriction on $t$ or $k$.

\subsection{An explicit height bound}
The geometric input is the Voronoi cell of real $E_8$.  Its vertices
have squared Euclidean norms $1$ and $8/9$; the first are deep holes
and twice each deep hole lies in $E_8$.
These facts follow from the two hole orbits, represented by
$(1,0,\ldots,0)$ and $(5,1,\ldots,1)/6$, respectively
\cite[Appendix~A.3.1]{bachoc2026cold}.  Thus the product of the two
horizontal Voronoi cells has $Q\leq1$.  Its vertices have $Q=1$
when both components are deep holes, and $Q\leq17/18$ otherwise.

\begin{theorem}[A cusp of height less than forty times the index]
\label{thm:ordinary-cusp-height}
Let $m\geq1$ and $v\in\Lambda$ satisfy $h(v,v)=3m$.
There is a primitive isotropic vector $e$ such that
\[
 0<N\!\left(\frac{h(e,v)}{\theta}\right)<40m.
\]
Equivalently, some unitarily equivalent vector has cusp coordinates
$v'=u e+\beta f+z$ with $0<N(\beta)<40m$.
This includes imprimitive positive vectors and vectors on root mirrors.
\end{theorem}
\begin{proof}
The height is a positive integer: lattice pairings are divisible by
$\theta$, and a nonzero isotropic vector cannot lie in the
negative-definite perpendicular of a positive vector.  Primitive
isotropic vectors exist, so choose $e$ minimizing this integer, and
write $D=N(\beta)$ in its cusp splitting.
Using \eqref{eq:cusp-inventory-translation}, reduce
$x=z/\beta$ into the product Voronoi cell.  Write
\[
 r=\frac mD,\qquad t=\frac u\beta=a+b\zeta,
 \qquad s=b+\frac a2\in[-1/2,1/2],
\]
where the last reduction uses a central translation.  The norm equation
gives $a=r+Q(x)$.  Minimality of the cusp $e$, compared with the primitive
cusp $f$, gives $N(t)\geq1$.  Since
\[
 N(t)=\frac34a^2+s^2,\qquad Q(x)\leq1,
\]
we obtain $1\leq a\leq1+r$ and $Q(x)\geq1-r$.

Suppose $r\leq1/40$.  Express $x=\sum_j\lambda_j y_j$ as a convex
combination of vertices of the product Voronoi cell.  The total weight
of vertices with $Q(y_j)<1$ is at most $18r$.  Moreover,
\[
 \sum_j\lambda_j Q(y_j-x)
   =\sum_j\lambda_j Q(y_j)-Q(x)\leq r.
\]
Some vertex $x_0$ with both components deep therefore satisfies
\[
 Q(x-x_0)\leq\frac r{1-18r}\leq2r,
\]
using $r\leq1/36$ for the last inequality.
Choose $s_0=1/2$ or $-1/2$ with the sign of $s$, and let
$t_0=1$ or $\bar\zeta$, respectively; these are the points with
$a=1$ and $s=s_0$.  The inequality $N(t)\geq1$ implies
\[
 \delta:=\frac12-|s|
 \leq3r+\frac32r^2,\qquad
 N(t-t_0)=\frac34(a-1)^2+\delta^2<16r^2.
\]
Twice every deep hole is integral.  Hence
\[
 w=2(t_0e+f+x_0)\in\Lambda
\]
is a nonzero isotropic vector, because $Q(x_0)=1$ and the $a$-coordinate
of $t_0$ is one.  With $v_0=t_0e+f+x_0$ we have $h(v_0,w)=0$.
The triangle inequality and positive Hermitian Cauchy--Schwarz give
\begin{align*}
 |h(v/\beta,w)|
 &\leq2\sqrt3\,|t-t_0|+2|h_+(x-x_0,x_0)|\\
 &\leq8\sqrt3\,r+6\sqrt{2r}\\
 &\leq\frac{\sqrt3}{5}+\frac3{\sqrt5}
 <\sqrt3=|h(v/\beta,e)|.
\end{align*}
If $w$ is imprimitive, dividing its Eisenstein content decreases its
pairing further.  This contradicts the choice of $e$.  Thus $r>1/40$
and $D<40m$.
\end{proof}

\subsection{Affine tests for the orbit representative}
\begin{proposition}\label{thm:pooled-affine-unranking}
In the cusp splitting there is a bounded system of representatives
with $H(\beta)^2<54m$, one for each ordinary open orbit of norm
$3m$. Apart from this height cutoff and the norm equation, membership
is a fixed Boolean combination of affine inequalities on the
twenty integral coordinates. The list of normals is independent of $m$.
\end{proposition}
\begin{proof}
Use the cusp splitting $v=(u,\beta,z)$ with $z\in K=F^2$ and
$H(a+b\zeta)=\max(|a|,|b|,|a+b|)$. We specify a presentation whose
tests remain affine as $\beta$ varies. Normalize by units and ordinary
conjugation to $\beta=a+b\zeta$, $a>0$, with either $b=0$ or
$a\leq b$. Here $H=H(\beta)=a+b$. In the integral Eisenstein basis
of the two caps, admit all $z$ whose sixteen real coefficients have
absolute value at most $H$. Every horizontal translation orbit has
such a representative: put each component of $z/\beta$ in the
centred coefficient parallelogram. Multiplication by $\beta$ then
bounds the two real coefficients by $(a+b)/2$ and $(a+2b)/2$,
respectively. Carry the Heisenberg scalar correction along with
this translation.

Write $u=p+q\zeta$. A central translation adds
$\beta\zeta=-b+(a+b)\zeta$. For $b>0$ choose $-H\leq 2q<H$;
for $b=0$ choose $-a\leq p+2q<a$. Each strip gives one point on
each integral central orbit, with the lower endpoint included. The norm equation
is $(a+b)p+bq=m+Q(z)$, so both strips give
$|u|=O(m/H+H)$ since $Q(z)=O(H^2)$. Also $Q(z/\beta)=O(1)$.
Order these presentations lexicographically by
$(H,\beta,z,u)$ in fixed real coordinates. Every period has positive
height at every primitive cusp. By \cref{thm:ordinary-cusp-height}, the bounds
$N(\beta_0)<40m$ and $N\leq H^2\leq 4N/3$ give a presentation
with $H^2<54m$. We may therefore impose the integer cutoff
$H\leq\lfloor\sqrt{54m-1}\rfloor$. There are finitely many admitted
tuples there, and the lexicographic minimum is unique in each orbit.

\emph{Finite comparisons throughout this box.}
For a null competitor $w=(A,\gamma,y)$ the transported coefficient is
$\beta'=h(w,v)/\bar\theta$. Suppose $H(\beta')\leq H(\beta)$
and write $D=N(\beta)$. Then $|h(w,v)|\leq2\sqrt D$.
The vector $w-(\gamma/\beta)v$ has zero $f$ coordinate, hence
its Hermitian norm is $-3Q(y-\gamma z/\beta)$. Expanding this norm
and bounding its cross term by $2|\gamma/\beta|\,|h(w,v)|$ gives
\[
 Q(y-\gamma z/\beta)\leq (4/3)|\gamma|-mN(\gamma)/D.
\]
Nonnegativity implies $|\gamma|\leq4D/(3m)<72$. A nonzero
$\gamma$ therefore requires $D\geq3m/4$.
Thus $y$ is bounded in the fixed positive lattice. If $\gamma=0$,
nullity forces $y=0$ and primitivity leaves only unit multiples of $e$.
Otherwise both $u/\beta$ and $z/\beta$ are bounded, and the equation for
$h(w,v)/\beta$ now bounds $A$. All bounds are effective from the
fixed Gram matrix and the displayed boxes. Enumerate this fixed
list once, retaining primitive null vectors. For each, compare
$H(\beta')$ and the unit- and conjugation-normalized coefficient
with the current key. These are fixed Boolean combinations of linear
signs, by splitting the coefficient plane into its twelve sectors.
A strictly smaller coefficient key always admits an actual presentation
after horizontal and central translation. Thus this test requires no
bound on the translations to a strictly smaller height.

For equality of the normalized coefficient, choose one fixed integral
transport of the competitor to the reference cusp. The single-cusp
orbit theorem supplies these transports, effectively by enumeration of integral isometries until each prescribed image occurs. Postcompose with a parabolic transformation
to compare the remaining coordinates. The parabolic group consists of
the Heisenberg translations, the finite cap group including exchange,
the units and the ordinary antiunitary map
$(u,\beta,z)\mapsto(\bar u,\zeta\bar\beta,c_K(z))$,
where $c_K$ is conjugation in a compatible integral cap frame.
Indeed an isometry fixing the primitive vector $e$ induces an
isometry on $e^\perp/\mathcal E e$. Removing its induced cap
action leaves a translation of $f$ by a cap vector, and after
\eqref{eq:cusp-inventory-translation} removes that vector, the
only freedom is a central translation. An isometry fixing the cusp
line may also multiply $e$ by a unit.
For a nonhorizontal tie, $D\asymp m$, so the input and its fixed
transport have all coordinates $O(H)$. In a competing presentation
the horizontal coordinates are again $O(H)$, which bounds the integral
translation $\tau$ because its displacement is $\beta\tau$ and
$|\beta|\asymp H$. Its scalar correction is $O(H)$. Since the target
strip has width one central period in its defining functional, only
boundedly many central shifts can give a competing presentation.

At the reference cusp this argument still bounds $\tau$, including
for $D\ll m$. The large transverse scalar in that case requires a
separate check. A unit preserving the chosen nonzero $\beta$ is the
identity. An antiunitary factor can preserve the preferred $\beta$
only on the two rays $b=a$ and $b=0$. On the diagonal it sends
$q$ to $-q$. On the real ray, after the normalizing unit it sends
$(p,q)$ to $(p,-p-q)$, reversing $p+2q$. The strip functional
therefore stays within $O(H)$ under all eligible finite factors;
the bounded horizontal shift changes it by $O(H)$ also. Its
increment in one central step has magnitude $2H$ in either strip.
This bounds the required central shifts uniformly even at thin cusps.
The transformations to test are consequently a fixed finite set.
Each is integral real linear in the \emph{whole} vector $v$:
in particular the Heisenberg correction for fixed $\tau$ is
$\beta Q(\tau)+h_+(\tau,z)/\theta$. A transformed tuple in the
box and strip with the same coefficient and a smaller $(z,u)$
is detected by affine signs. Identical lexicographic keys from
stabilizers retain the same integer tuple once.

For completeness, mirror exclusion is affine here too. For a root
$w=(A,\gamma,y)$ of norm $-3$ perpendicular to $v$, subtract
$(\gamma/\beta)v$ to obtain
$Q(y-\gamma z/\beta)=1-mN(\gamma)/D$. This bounds $\gamma,y$.
If $\gamma\ne0$ then $D\geq m$, so $h(w,v)=0$ bounds $A$ by the
same scalar estimate. If $\gamma=0$ it gives $Q(y)=1$ and bounds
$A$ directly from bounded $z/\beta$. Enumerate those roots and
exclude their perpendiculars. The ordinary open period correspondence
now gives exactly one surviving tuple for each requested graph,
including imprimitive periods.
\end{proof}

For a retained tuple define $s=q$ on the branch $b>0$ and
$s=p+2q$ on the branch $b=0$. Its code is
\[
 (H,a,z_1,\ldots,z_{16},s),\qquad H=a+b.
\]
The arithmetic rank of a cage is the position of its retained code in lexicographic order, starting at zero. A coefficient family consists of the cages of the given size with the same pair $(a,b)$. We put
$(H,a)$ first so that all isomers with the same coefficient $\beta$
occupy one interval. Its endpoints can be computed once and reused
across rank queries at the same size: subtracting the lower endpoint
gives the rank within that coefficient
(\cref{cor:orthoscheme-height-directory}). At each
height the skew coefficients occur in increasing $a$, followed by
the real coefficient $(H,0)$. Within one coefficient compare the
first cap and then the second in the fixed integral basis. For each
cap pair visit the surviving scalar values in increasing $s$. This
also orders them by increasing $q$ on the real branch, since there
$p=(m+Q(z))/H$ is fixed. The code recovers $b=H-a$ and the norm
equation recovers $p$ and $q$. Thus each ordinary isomer has exactly
one rank. The earlier minimum in $(H,\beta,z,u)$ chooses its
representative; this code orders the representatives so chosen.

\begin{figure}[!t]
\centering
\includegraphics[width=.76\linewidth]{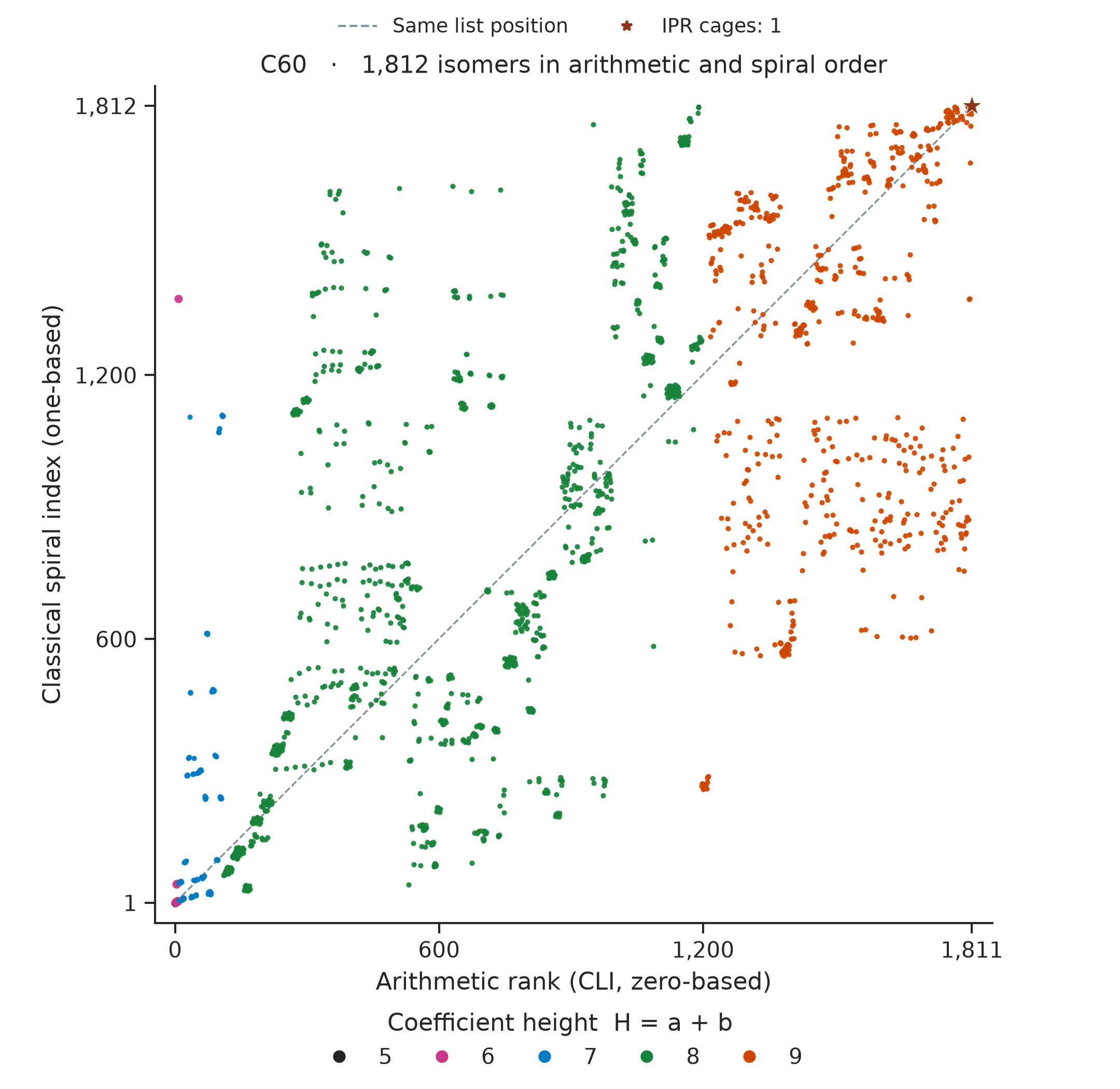}
\caption{One point for each of the $1{,}812$ ordinary $\mathrm C_{60}$
isomers. The horizontal coordinate is its zero-based arithmetic rank,
the vertical coordinate its one-based classical spiral index, and
the color gives $H=a+b$. The dashed line is $y=x+1$. The star is the
unique cage with isolated pentagons, last in both orders.}
\label{fig:c60-spiral-ranks}
\end{figure}

For comparison, the classical spiral index of a $\mathrm C_{60}$
cage is its one-based position in the ordering by pentagon-position
tuples \cite[Sec.~3.3]{bille2025open}. For each cage take the lexicographically first tuple of
twelve pentagon positions among the no-jump spirals of each cage,
allowing either orientation and any starting face
\cite[Sec.~III]{wirz2018faceSpirals}. Sorting these tuples gives the
vertical coordinate in \cref{fig:c60-spiral-ranks}. The
arithmetic order groups each height into a single interval; a spiral
index can separate cages of the same height.

\section{Developing the graph from its period}
\subsection{A single star polygon in twenty integer coordinates}
\label{sec:constructive-star-domain}

Thurston cuts from one cone vertex to the other eleven and
develops the resulting disk \cite[Section~7 and Proposition~7.1]{thurston1998}.
Its boundary gives a fixed-size integral chart.

Put $\mathcal E=\mathbb Z[\zeta]$, where
$\zeta=e^{\pi i/3}$, and write
$N(a+b\zeta)=a^2+ab+b^2$ and
$[a+b\zeta,c+d\zeta]=ad-bc$. Choose
\[
 p_0=0,\qquad p_1,\ldots,p_{10}\in\mathcal E,
 \qquad p_{11}=p_0,
\]
and define, with subscripts modulo eleven,
\[
 e_i=p_{i+1}-p_i,\qquad q_i=p_i+\zeta e_i,
 \qquad Q=(p_0,q_0,p_1,q_1,\ldots,p_{10},q_{10}).
\]
Require $Q$ to be a simple counterclockwise polygon and impose the
finite list of quadratic inequalities
\begin{equation}\label{eq:star-nearest-tip-domain}
 N(q_i-p_j)\geq N(e_i)\qquad(0\leq i,j\leq10).
\end{equation}
The two comparisons with $j=i,i+1$ are equalities automatically.
Simplicity excludes zero sides and unintended boundary contacts.

\begin{proposition}[Integral star construction]
\label{prop:integral-star-construction}
For every such array, identify $[p_i,q_i]$ with
$[p_{i+1},q_i]$, matching the points at equal distance from $q_i$.
The quotient is a Euclidean cone sphere with twelve distinct cone
points, all of angle $5\pi/3$. Its area in elementary triangular-lattice
triangles is the positive even integer
\begin{equation}\label{eq:star-area-shell}
 V=\sum_{i=0}^{10}[p_i,p_{i+1}]
       -\sum_{i=0}^{10}N(p_{i+1}-p_i).
\end{equation}
The paired sides descend to shortest paths from the common image of
the $p_i$ to the eleven points $q_i$. Its unit mesh is the dual of
a classical fullerene on $V$ vertices. Every classical fullerene dual
has at least one representation satisfying these conditions.
\end{proposition}
\begin{proof}
Since $1-\zeta=\zeta^{-1}$, the turn at $q_i$ along the boundary
of $Q$ is $-2\pi/3$. Its interior angle is therefore $5\pi/3$.
All $p_i$ become one vertex, while the eleven $q_i$ stay distinct.
Their quotient links are circles, and the oriented connected quotient
has twelve vertices, eleven edges and one face, hence is a sphere.
The sum of its angles at the common root is
\[
 (22-2)\pi-11\cdot5\pi/3=5\pi/3.
\]
Subtracting the eleven equilateral notches from the signed area of
the polygon with vertices $p_i$ gives \eqref{eq:star-area-shell}.
This is a signed-area identity and does not require that latter
polygon to be convex. Modulo two, the terms linear in the coordinate
differences telescope, while the mixed terms in the sum of norms
equal the sum of determinants. Thus $V$ is even. Every side pairing
is an Eisenstein translation followed by a lattice unit, so the unit
triangular tessellation is compatible with all identifications.

For the shortest-path assertion, let
$f(x)=\min_j|x-p_j|$ on $Q$. The Euclidean Voronoi cell of $p_i$
is convex. By \eqref{eq:star-nearest-tip-domain}, it contains both
$p_i$ and $q_i$, hence their segment. The corresponding assertion
holds for $p_{i+1}$. Consequently $f$ has equal values at paired
boundary points. It descends to a one-Lipschitz function on the
quotient, vanishing at the root. At $q_i$ it equals $|e_i|$, and
the paired side is a path of exactly that length.

Conversely, choose a cone vertex of a fullerene dual as root and
shortest paths from it to the remaining cone vertices. A shortest
path cannot pass through a positive-curvature cone in its interior.
Two such paths cannot cross away from the root: switching their
prefixes at an intersection and smoothing a resulting corner would
shorten one. Cutting along them gives Thurston's star disk. His
Proposition~7.1 embeds its enlarged flower in the plane, and hence
embeds the disk as a simple polygon. Development preserves the
Eisenstein lattice, and the cone angles give precisely the formula
for $q_i$ above.

For completeness of the nearest-tip condition, attach to each slit
pair the circular sector of angle $\pi/3$ and radius equal to the
slit's length, as in Thurston's flower. Map the sector radially onto
the original slit and map the star disk back to the sphere. These
maps agree on their interfaces and do not increase lengths; the
entire flower boundary maps to the root. The intrinsic distance
from $q_i$ to that boundary is therefore at least its shortest-path
distance $|e_i|$ to the root. Since the flower is embedded in the
plane, the open Euclidean disk of that radius about $q_i$ is contained
in it: otherwise the first exit of a straight segment would reach
the boundary at a smaller intrinsic distance. All the $p_j$ are on
that boundary, proving \eqref{eq:star-nearest-tip-domain} even when
a segment from $q_i$ to another tip leaves the original star disk.
The mesh and simplicity assertions follow from the next lemma.
\end{proof}

\begin{lemma}[The integral cone mesh is classical]
\label{lem:star-classical-mesh}
A compact Euclidean cone sphere with twelve cone angles $5\pi/3$,
Eisenstein affine holonomy, and all cone points on the developed
lattice has a simple unit triangular mesh, whose dual is a classical
fullerene.
\end{lemma}
\begin{proof}
The affine holonomy preserves the triangular tessellation, so its
pullback is intrinsic. Cone points lie at mesh vertices. Away from
the grid, local development continues along every path in an open
unit triangle: a lift cannot terminate at a cone or a grid edge
while the path stays a positive distance from the triangle boundary,
and compactness supplies continuation. Thus the universal cover of
an open-cell component develops onto an open unit triangle. The
triangle is simply connected, so this is an isometry. Any nontrivial
remaining cell identification would be an orientation-preserving
triangle rotation and would fix its centre, creating a further cone
point. There is none. Hence open cells are ordinary unit triangles;
their closures may initially have boundary identifications. Finite
area gives finitely many cells. The resulting spherical triangular
map has twelve degree-five vertices and all other vertices of
degree six.

To prove simplicity, suppose there is a loop or a pair of parallel
edges. On each side of that separating curve, choose a disk bounded
by a loop or digon with the fewest triangular faces. Put $b=1$ or
$2$ for its boundary length, $I$ for its number of interior vertices,
$k$ for its number of interior degree-five vertices, and $D$ for
the sum of its boundary degrees within the disk. Minimality excludes
all further loops and parallel edges away from its boundary: such
edges would enclose a smaller eligible disk. Every interior vertex
therefore has at most $I-1+b$ neighbours. Such a vertex exists,
since the number of faces is $2I+b-2>0$. Its degree is at least five, giving
$I\geq6-b\geq4$. Disk Euler characteristic and the triangle
incidences give
\[
 k+4b-D=6,\qquad
 E_{\rm int}=\frac{6I-k-D+2b}{2},
\]
where $E_{\rm int}$ counts edges between interior vertices.
That induced graph is simple and planar, so
$E_{\rm int}\leq3I-6$. It follows that
$k+D-2b\geq12$, and hence $k\geq9-b\geq7$.
The two chosen disks have disjoint interiors and would require at
least fourteen degree-five vertices, a contradiction. A simple
spherical triangulation is a maximal planar graph and is
three-connected; dualization now gives the claimed fullerene.
\end{proof}

\begin{proposition}[The star chart has the common period lattice]
\label{prop:star-marked-period-chart}
An integral marking of
the star chart may be chosen as
\begin{equation}\label{eq:star-ambient-marking}
 \Phi(p)=\sum_{j=0}^{9}p_{10-j}\alpha_j.
\end{equation}
It satisfies $h(\Phi(p),\Phi(p))=3V/2$. Consequently an ambient
vector whose reversed coordinates satisfy the star-domain conditions
has an explicit marked geometric return.
\end{proposition}
\begin{proof}
The outgoing slit vectors are $w_i=q_i-p_i=\zeta e_i$.
Their paired boundary vectors are $\zeta^{-2}w_i$, and closure of
the developed boundary is equivalent to $\sum_iw_i=0$.
Thus the formal twisted period module of this cut diagram is
$\ker(\sum:\mathcal E^{11}\to\mathcal E)$. Its internal diagonals
are determined by differences of developed vertex positions; they
add no relation or extra parameter. The change
\[
 p_i=\zeta^{-1}\sum_{j=0}^{i-1}w_j
\]
between its ten free slit vectors and $(p_1,\ldots,p_{10})$ is
integral with integral inverse. Hence these are coordinates for the
whole period lattice of the chart, rather than a sublattice obtained
by clearing denominators. This is also the leader-and-follower
construction of \cite[Proposition~3.2]{thurston1998}: the final
holonomy requires division only by the unit $1-\zeta^{\pm1}$.

Expanding \eqref{eq:star-area-shell}, the Hermitian form $3V/2$
on the ten coordinates $p_i$ has diagonal $-3$, upper adjacent
entry $1+\bar\zeta=2-\zeta$, lower adjacent entry $1+\zeta$,
and zero other entries. Reversing their order gives exactly the
Gram matrix of $\Lambda$. This identifies the entire
integral period module and its form, not only the norm of a single
candidate. Restricting this marking to the open set of nondegenerate
star diagrams gives \eqref{eq:star-ambient-marking} for their actual
geometric periods.
\end{proof}

For a prescribed developing frame the scalar unit in this identification
matters. Write $g_i(z)=q_i+\zeta(z-q_i)$, so $g_i(p_i)=p_{i+1}$.
The composite $g_{10}\cdots g_0$ fixes $p_0$ with linear part
$\bar\zeta$. The rotation about $p_0$ completes its product to the
identity. A fan of paths inside the cut disk therefore reads the ordered
points $Q_j=q_{10-j}$ for $0\leq j\leq10$ and $Q_{11}=p_0$, in the
rightmost-first convention of \cref{cnt-sec:integral-representatives} of~\cite{fullereneCountManuscript}.
Put $\theta=1+\zeta$. The ten chain pairings give
\[
 h(\alpha_j,\bar\zeta\Phi(p))
    =\theta(q_{9-j}-q_{10-j})
    =\theta(Q_{j+1}-Q_j),\qquad 0\leq j\leq9 .
\]
Thus this frame reads the literal period $\bar\zeta\Phi(p)$. Rotating
the frame by $\zeta$ gives the convention of
\cref{prop:star-marked-period-chart}. Once the final star is developed,
reversing its recorded changes of marking returns the input vector
in its original frame.

\subsection{Comparisons of slit lengths}
Set $E(p)=\sum_{i=0}^{10}N(p_{i+1}-p_i)$.
The adjacent Hurwitz move
\begin{equation}\label{eq:star-hurwitz-coordinate-move}
 \begin{aligned}
 (g_i,g_{i+1})&\longmapsto
 (g_{i+1},g_{i+1}g_i g_{i+1}^{-1}),\\
 p_{i+1}&\longmapsto p_{i+2}+\zeta(p_i-p_{i+1})
 \end{aligned}
\end{equation}
leaves all other tips unchanged.  Its integral linear matrix has order
three and preserves the star area form.  In root-chain coordinates it
is an integral unitary triflection.

\begin{proposition}[Descent to the weak nearest-tip inequalities]
\label{prop:star-hurwitz-energy}
If $N(q_i-p_j)<N(q_i-p_i)$, a sequence of adjacent Hurwitz moves or
their inverses replaces the $i$th slit length by $|q_i-p_j|$ and
preserves the other slit lengths, up to reordering.  Consequently,
the total integer energy strictly decreases.  Repeating such moves
terminates at the weak inequalities
\[
 N(q_i-p_j)\geq N(p_{i+1}-p_i)\qquad(0\leq i,j\leq10).
\]
This termination statement does not assert that the resulting polygon
is simple.
\end{proposition}
\begin{proof}
Under \eqref{eq:star-hurwitz-coordinate-move}, the new centers are
$q_{i+1}$ and $g_{i+1}(q_i)$, and the middle tip becomes
$g_{i+1}(p_i)$.  The old $i$th slit retains its length because
$g_{i+1}$ is a rotation.  The other slit changes from
$|q_{i+1}-p_{i+1}|$ to $|q_{i+1}-p_i|$.  Moving a selected factor
left across successive earlier factors therefore gives every earlier
root-copy comparison.  The inverse move similarly gives every later
comparison; the adjacent endpoint comparisons are equalities.
There are at most ten adjacent moves in each such change.
Every accepted change lowers the nonnegative integer $E$, so the
process terminates after at most its initial value many changes.
\end{proof}

\subsection{The geometry of a weak minimum}
\label{sec:star-weak-minimum-geometry}
Here a \emph{weak minimum} means an array satisfying all nearest-tip
inequalities, not merely one on which the twenty single adjacent moves
fail to lower $E$.  A comparison change in
\cref{prop:star-hurwitz-energy} can pass through intermediate arrays
of larger energy.  Only its endpoint energy is used in the descent
argument.

\begin{lemma}[Coincidences excluded by mirror-freeness]
\label{lem:star-coincidences}
In a positive mirror-free star period, every slit has nonzero length
and the eleven points $q_i$ are pairwise distinct.  If two root tips
coincide, their shorter cyclic separation is two or three.
\end{lemma}
\begin{proof}
Translation followed by cyclic reindexing of the tips preserves the
area form and the Eisenstein lattice, so it is an integral unitary
change of marking.  A coincidence can therefore be written
$p_0=p_k=0$, where $1\leq k\leq5$.
In ordered pairs of real coordinates on the tips, the symmetric
matrix $B$ representing $V$ has diagonal blocks
$\bigl(\begin{smallmatrix}-2&-1\\-1&-2\end{smallmatrix}\bigr)$ and adjacent
upper blocks $\bigl(\begin{smallmatrix}1&1\\0&1\end{smallmatrix}\bigr)$.
The operator $J$ is multiplication by $\zeta$ on each coordinate.
Block tridiagonal inversion gives the
following diagonal blocks of $B^{-1}$, in order for $p_1$ through
$p_5$:
\[
 M,\quad 2M,\quad 2M,\quad M,\quad 0,
 \qquad M=\begin{pmatrix}-2&1\\1&-2\end{pmatrix}.
\]
These identities also follow by the tridiagonal determinant recurrence
$D_n=-3D_{n-1}-3D_{n-2}$ for the Hermitian form, with
$D_0=1$ and $D_1=-3$.  Let $\epsilon$ select the first real
coordinate of $p_k$ and put $r=B^{-1}\epsilon$.  The determinant is one (equivalently use
$|\det H_\star|=3^5$ and the real Gram determinant
$3^{10}|\det H_\star/3|^2=1$, with $H_\star$ the Hermitian matrix
in tip order). Unimodularity makes
$r$ integral.  If $p_k=0$, then both $B(r,v)$ and $B(r,Jv)$ vanish.
For $k=1,4$, the displayed block gives $B(r,r)=-2$, so $r$ is a
perpendicular root.  For $k=5$, it gives a nonzero isotropic vector
in the Hermitian perpendicular of a positive vector, contradicting
negative definiteness.  This excludes separations one, four and five,
and in particular excludes zero-length slits.

Suppose $q_i=q_j$ with $i<j$.  Move the $j$th factor left until
these two factors are adjacent.  Its centre stays fixed; the old
$i$th centre is rotated about that same point and therefore stays
there too.  In star coordinates, adjacent centres satisfy
\[
 h_\star(e_k,v)
  =(1-2\zeta)(q_k-q_{k-1}),
\]
where $e_k$ is the $k$th coordinate vector, of norm $-3$.
Their equality therefore gives a root mirror.  Hurwitz moves are
integral unitary maps, so mirror-freeness excludes the original
coincidence as well.
\end{proof}

\begin{proposition}[Only two- and three-slit bubbles can obstruct simplicity]
\label{prop:star-small-bubbles}
Let $v$ be a positive mirror-free weak minimum.  The drawn boundary
$Q$ is an embedded planar Eulerian network after equal root tips are
identified.  If $Q$ is not a simple polygon, it has a cyclically
contiguous subwalk
\[
 P,p_{i+1},\ldots,p_{i+k-1},P
\]
of root tips with $k=2$ or $3$, whose alternating tip-and-valley
boundary is a clockwise Jordan polygon.  We call this subwalk a
negative bubble.

More generally, decomposing the directed network into Jordan cycles
gives exactly one counterclockwise cycle and at most four slit edges
in clockwise cycles.  This assertion does not require the network
to be a cactus or assume that all its cycles can be reordered by
energy-preserving moves.
\end{proposition}
\begin{proof}
Use the ordinary Euclidean Voronoi cells of the set of distinct tip
positions.  Each segment from a tip to an incident $q_i$ stays in
its cell; its open segment lies in the cell interior, since the tip
itself is strictly closer to that site than to any other distinct
site.  Segments from different sites therefore have disjoint
interiors.  Two segments from the same site can overlap only if
they end at the same boundary point of its convex Voronoi cell:
a ray from an interior point meets that boundary at most once.
\Cref{lem:star-coincidences} excludes this.  Each $q_i$ has degree
two, while every distinct tip site has equal incoming and outgoing
degrees.  Thus the only remaining self-contacts of $Q$ are repeated
tip sites in a planar Eulerian network.

A directed Jordan cycle with $k$ valleys has the prescribed turn
$-2\pi/3$ at each of them.  If its orientation is counterclockwise,
the sum of its tip angles is
\[
 (2k-2)\pi-k\cdot5\pi/3=(k/3-2)\pi>0.
\]
It must have $k\geq7$.  Signed area is additive under decomposition
of the directed Eulerian walk into Jordan cycles.  Positive total
area requires a counterclockwise cycle, and eleven available slit
edges permit only one.  At most four slit edges remain in clockwise
cycles.

If the boundary is not simple, choose equal tips at shorter cyclic
distance two or three, as provided by
\cref{lem:star-coincidences}.  Its interior tip occurrences are
distinct: a repetition would give a still shorter coincidence,
including a zero-length edge in the three-slit case.  The resulting
embedded closed subwalk is a Jordan polygon.  Since it has fewer
than seven valleys, it is clockwise.  This identifies a contiguous
bubble directly, without a reordering argument.
\end{proof}

\subsection{A bounded escape from every two-slit bubble}
\label{sec:star-two-bubble-escape}

\begin{proposition}[Two zero-energy moves suffice]
\label{prop:star-two-bubble-escape}
A positive mirror-free weak minimum containing a contiguous
two-slit bubble admits at most two adjacent, energy-preserving
Hurwitz moves after which some nearest-tip comparison is strict
in the wrong direction.  Applying that comparison change then
strictly lowers $E$.  The comparison change may have intermediate
states of larger energy; its endpoint has smaller energy.
\end{proposition}
\begin{proof}
Cyclically reindex and translate so the bubble tips are
$0,x,0$, with $x\ne0$.  Let the next tip outside the bubble be
$y\ne0$.  The internal order-three move rotates $x$ through
\[
 x,\quad \zeta^2x,\quad \zeta^4x
\]
in one of the two cyclic orders.  The two bubble slit lengths
remain $|x|$ and all other slit lengths remain unchanged.
Its valleys are $\zeta x$ and $\zeta^{-1}x$; the following
external valley is $\zeta y$.
\Cref{fig:two-slit-recut} shows how a rotation exposes a closer tip.

\begin{figure}[htbp]
\centering
\includegraphics[width=.88\linewidth]{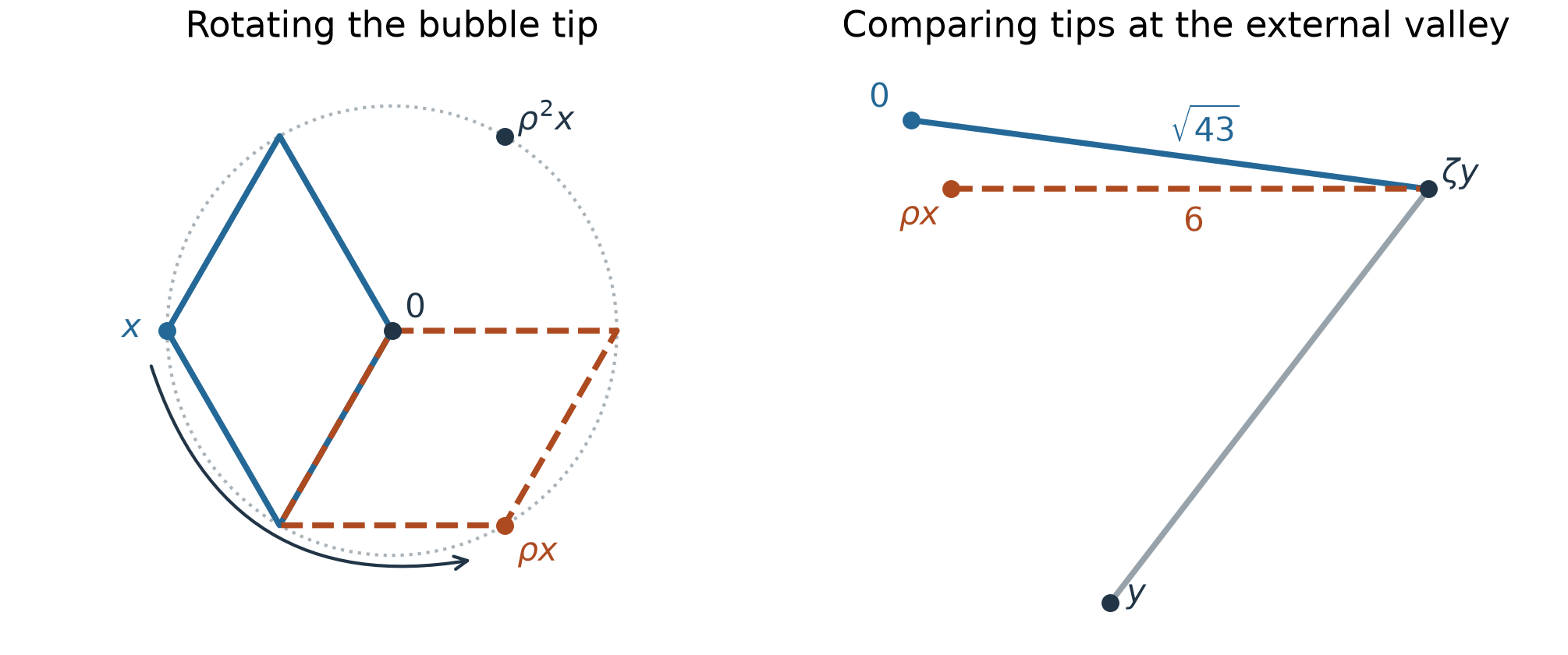}
\caption{One two-slit recut from the $232$-triangle example in
\protect\path{data/star_recuts.json}. In unit Eisenstein coordinates
the initial tips are $0,x,0$ with $x=-1$, and the next tip is
$y=6-7\zeta$. Rotating $x$ by $\rho=\zeta^2$ preserves the two
slit lengths. The external valley $\zeta y$ now has squared distance
$36$ to $\rho x$, compared with $43$ to its assigned tip $0$.
The rotation preserves energy; the ensuing comparison lowers it.}
\label{fig:two-slit-recut}
\end{figure}

Suppose all three rotated arrays still satisfy the weak inequalities.
Write $\rho=\zeta^2$ and
\[
 M_0=\max_{0\leq j\leq2}
       \operatorname{Re}\bigl(\zeta y\,\overline{\rho^jx}\bigr).
\]
Comparing the external valley against every rotated bubble tip gives
$N(x)\geq2M_0$.  Comparing every rotated bubble valley against
the external tip gives $N(y)\geq2M_0$: the three odd powers of
$\zeta$ in these comparisons are the same set.
Three directions spaced by $2\pi/3$ have maximal cosine at least
$1/2$, so $2M_0\geq|x||y|$.  The two preceding inequalities force
$|x|=|y|$ and equality throughout.  Equality in the cosine bound
then gives
\[
 y/x\in\{1,\zeta^2,\zeta^4\}.
\]
In one rotated array the bubble tip is therefore $y$.  Its valley
$\zeta y$ coincides with the external valley, contradicting
\cref{lem:star-coincidences}.  At least one of the two other
rotations must violate a nearest-tip inequality.  Reach the first
such rotation, then use \cref{prop:star-hurwitz-energy}.
\end{proof}

\subsection{Six recuts repair the three-slit case}
\label{sec:star-three-bubble-escape}

\begin{lemma}[A three-slit bubble has no other attachment]
\label{lem:star-three-bubble-contacts}
Suppose a positive mirror-free weak minimum has a contiguous
three-slit bubble and has no contiguous two-slit bubble.  The
negative bubble and its complementary eight-slit boundary are
Jordan polygons meeting only at their common initial tip.  The
complementary polygon is counterclockwise and has larger unsigned
area than the bubble.
\end{lemma}
\begin{proof}
Normalize the bubble tips to $p_0=p_3=0$, with $a=p_1$ and
$b=p_2$.  Its complement has eight slit edges and positive signed
area.  Decomposing that complement into directed Jordan cycles
requires a positive cycle with at least seven edges.  A second
cycle would require at least two more, so the complement is a single
counterclockwise Jordan polygon.  Positive total area says its area
exceeds the bubble's unsigned area.

By \cref{lem:star-coincidences}, coincidences of $a$ with exterior
tips can occur only at $p_4,p_9,p_{10}$, and coincidences of $b$
only at $p_4,p_5,p_{10}$.  The equalities $a=p_{10}$ and $b=p_4$
would be two-slit bubbles.  The equalities $a=p_4$ and $b=p_{10}$
would give respectively $q_0=q_3$ and $q_2=q_{10}$, which are
root mirrors.  It remains to exclude $a=p_9$ and $b=p_5$.

Let $H_\star$ be the Hermitian matrix in star coordinates, with
diagonal $-3$ and upper adjacent entry $\bar\theta$, where
$\theta=1+\zeta$.  Let $e_k$ select $p_k$, and set
\[
 A=H_\star^{-1}e_3,\qquad
 D=H_\star^{-1}(e_9-e_1),\qquad
 F=H_\star^{-1}(e_5-e_2).
\]
The tridiagonal inverse gives
\[
 h(A,A)=h(D,D)=h(F,F)=-2,\qquad
 h(A,D)=-\bar\theta,\qquad h(A,F)=\theta.
\]
Thus $h(A-D,A-D)=h(A+F,A+F)=-1$.
The matrix $H_\star/\theta$ is integral with unit determinant:
$\det H_\star=-3^5$ and $\theta^{10}=3^5\zeta^5$.
Consequently $\theta\Lambda^\vee=\Lambda$, and
\[
 r_1=\theta(A-D),\qquad r_2=\theta(A+F)
\]
are integral roots of norm $-3$.  The conditions $p_3=0$ and
$p_9-p_1=0$ make $r_1$ perpendicular to the period.  The conditions
$p_3=0$ and $p_5-p_2=0$ do the same for $r_2$.  Mirror-freeness
excludes both contacts.  Distinct arms cannot intersect away from
tip sites by the Voronoi argument, so only the common initial tip
remains.
\end{proof}

\begin{proposition}[Bounded escape from a three-slit bubble]
\label{prop:star-three-bubble-escape}
A positive mirror-free weak minimum with a contiguous three-slit
bubble admits at most five applications of a fixed two-move
energy-preserving recut after which either a nearest-tip inequality
fails or a two-slit bubble is present.  In the latter case
\cref{prop:star-two-bubble-escape} applies.  Thus a bounded number
of recuts always leads to a strict decrease of the endpoint energy.
\end{proposition}
\begin{proof}
Again write the bubble tips as $0,a,b,0$, and put $\rho=\zeta^2$.
The inverse adjacent moves at its two internal positions give
\begin{equation}\label{eq:star-three-bubble-cycle}
 C(a,b)=\rho(a-b,a),\qquad C^3(a,b)=(-a,-b),\qquad C^6=1.
\end{equation}
The three slit lengths are permuted, so each whole application of
$C$ preserves $E$.  Its intermediate adjacent move need not do so.
Geometrically it sends the negative hexagon to a rigid copy:
\[
 z\longmapsto\zeta^{-1}(z-a).
\]
The old tip $a$ becomes its initial tip, and the old cyclic sequence
$(a,b,0)$ becomes
$(0,\zeta^{-1}(b-a),-\zeta^{-1}a)=(0,\rho(a-b),\rho a)$.
The valley coordinates transform by the same rigid motion.  The
six states therefore include each of the three possible initial
tips, together with its opposite copy about the common initial tip.

Suppose every one of these six states is weak and contains no
two-slit bubble.  The three valley angles of the clockwise hexagon
are $\pi/3$, so its three tip angles sum to $3\pi$.  One tip angle
is at least $\pi$.  Choose the corresponding state and its opposite.
Their open interior angular sectors at the initial tip cover every
direction, except possibly their boundary rays when the angle is
exactly $\pi$.  An external arm cannot lie on a bubble boundary
ray: both endpoints would be on the same ray from a Voronoi site
to its cell boundary, forcing coincident valleys.

In at least one of the two states an external arm therefore enters
the interior of the negative bubble.  By
\cref{lem:star-three-bubble-contacts}, the complementary positive
Jordan polygon has no other contact with the bubble boundary.
It cannot cross that boundary, by the Voronoi embedding.  Its entire
boundary, and hence its enclosed region, must consequently lie
inside the negative bubble.  This contradicts the strict area
inequality in that lemma.  At least one of the other five states
must instead violate a nearest-tip inequality or contain a two-slit
bubble, as asserted.
\end{proof}

\subsection{A terminating recut normalization}
\label{sec:star-complete-recut-normalization}

\begin{theorem}[Reconstructing a marked star polygon from its period]
\label{thm:star-complete-recut-normalization}
Every positive mirror-free integral period can be transformed into
an admissible simple star polygon by a deterministic sequence of
cyclic changes of starting tip and Hurwitz recuts.  If $E_0$ is its
initial slit energy, the normalization uses $O(E_0)$ fixed-size
integer operations, with $O(\log(2E_0))$-bit coordinates, and records
a transformation word of length $O(E_0)$.  Its bit cost is
$O(E_0(\log(2E_0))^c)$ for a fixed $c$.

For $h(v,v)=3m$, expanding the normalized star gives the associated
classical fullerene, with an additional $O(m\log m)$ word cost.
The recorded inverse word supplies a marking whose period is the
original $v$.  The theorem does not require an independently
generated graph census or a search through unitary words.
\end{theorem}
\begin{proof}
Fix orders on tip pairs and on the finitely many recuts considered
below.  Apply the first violated nearest-tip comparison until a
weak minimum is reached.  If it is simple, stop.  Otherwise
\cref{prop:star-small-bubbles} gives a contiguous two- or three-slit
bubble.  Cyclically normalize its initial tip.  In the two-slit case,
use \cref{prop:star-two-bubble-escape}.  In the three-slit case,
inspect the six states of \eqref{eq:star-three-bubble-cycle} in the
fixed order, stopping at a violated comparison or a two-slit bubble;
the latter has the already proved escape.  Follow the violated
comparison to an endpoint of strictly smaller energy and repeat.
All selected recuts preserve the integral lattice, the area and
mirror-freeness.  The energy is a nonnegative integer, so at most
$E_0$ strict endpoint decreases are possible.

There is a fixed bound on the work between two such decreases:
eleven tips, $121$ nearest-tip comparisons, six three-bubble states,
three two-bubble states, and at most ten signed adjacent moves in
a comparison change.  No branch of unbounded depth is explored.
At an endpoint with $p_0=0$, Cauchy--Schwarz gives
\[
 |p_i|\leq\sum_{j<i}|p_{j+1}-p_j|\leq\sqrt{11E_0}.
\]
The bounded number of intervening integral linear recuts increases
this bound by only a fixed factor.  This proves the stated arithmetic
and bit bounds.  The transformation is retained as a word of bounded
length per decrease; the argument does not assume small entries in
its fully multiplied unitary matrix.

When the algorithm stops, all tips are distinct, the drawn boundary
is simple by the Voronoi argument, and its positive area makes its
orientation counterclockwise.  It is therefore an admissible star
polygon.  The integral star construction and classical-mesh lemma
(\cref{prop:integral-star-construction,lem:star-classical-mesh})
give its graph. Each slit is shortest, and its endpoints are mesh
vertices joined by a simple edge path of length at most the number
of mesh vertices minus one, hence $O(m)$. The polygon therefore has
$O(m)$ perimeter and coordinates. The strip construction of
\cref{cor:linear-census-adjacency} gives the stated expansion bound.  If $g$ is the recorded transformation,
the resulting star has coordinate period $gv$; assigning it the
inverse transported marking restores period $v$ exactly.
\end{proof}

\section{Cost of reconstruction}\label{sec:decoder}
\begin{theorem}[Linear recut delay for every compressed member]
\label{thm:linear-recut-delay}
Let $v=ue+\beta f+z$ be a centered ordinary open period of positive
index $m$ in the cusp splitting, with $D=N(\beta)\le C_0m$ for a fixed
$C_0$. The bounded-recut decoder reaches an admissible star in $O(m)$
recut letters and $O(m)$ arithmetic word operations, uniformly in the
period, its phase and its height. The words hold $O(\log(2m))$ bits.
\end{theorem}
\begin{proof}
On anchored tips $p_0=0,p_1,\ldots,p_{10}$ let
$E=\sum_i|p_{i+1}-p_i|^2$ with cyclic subscripts. Its square root
is a norm on the period space. We first claim that the reference cusp
$e$ realizes the absolute minimum of $E$ on nonzero integral null
vectors. For a null vector the signed axial shoelace sum equals $E$ by
the star-area identity. In the ordinary Euclidean plane the algebraic
area is therefore $\sqrt3 E/4$. Every nonzero integral edge has
length at least one, so the perimeter $P$ is at most $E$.
The signed isoperimetric inequality gives $\sqrt3 E/4\le P^2/(4\pi)$
for this possibly self-intersecting loop: parametrize by arclength,
subtract its mean, and apply Cauchy--Schwarz and the periodic
Wirtinger inequality to the signed-area integral. Hence the positive
integer $E$ is at least $\pi\sqrt3>5$. Conversely take six consecutive unit edges of a counterclockwise
regular hexagon and five zero edges. Its shoelace sum and energy
are both six. It gives a primitive integral null vector (one tip
is a unit). By the single-cusp theorem we may and do use this vector
as the reference $e$ in the cusp-to-star identification. An integral
isometry transporting the fixed splitting to this cusp is independent
of $m$ and can be chosen once by enumeration.

The box gives $Q(z)=O(D)$. Write $u/\beta=A+B\zeta$ with real
$A,B$. The norm equation is $A=(m+Q(z))/D$ and central translation
puts $B+A/2$ in $[-1/2,1/2)$. We call this a centered representative.
Since $N(A+B\zeta)=3A^2/4+(B+A/2)^2$ and $D=O(m)$, we obtain
$|u|^2=3m^2/(4D)+O(m)$ and $\|\beta f+z\|_E=O(\sqrt D)$.
Expanding the fixed Hermitian energy gives
\[
 E(v)=\tfrac{9}{2}m^2/D+O(m).
\]
Let $w=Gv$ be the final admissible star produced by descent and let
$n=Ge$. Integral recuts and cyclic changes preserve the Hermitian
pairing and $n$ is nonzero and null. We have $E(w)\le E(v)$ and
$|h(w,n)|=\sqrt{3D}$. Only a fixed finite set of $n$ is possible.
Indeed, in a fixed complex basis diagonalizing $h$ with signs $(1,9)$,
write $w=(w_0,w_-)$ and $n=(n_0,n_-)$. Then
\[
 |h(w,n)|\ge |n_0|(|w_0|-\|w_-\|)
 =\frac{3m|n_0|}{|w_0|+\|w_-\|}.
\]
Here $|n_0|=\|n_-\|$ and the denominator is $O(\sqrt{E(w)})$.
Thus $\|n\|=O(\sqrt{D E(v)}/m)=O(1)$, and integrality proves finiteness.

For each such image choose one integral transport $T_n$ of $e$ to $n$.
The finite set of frames gives uniform norm-comparison constants. Write
$T_n^{-1}w=u'e+\beta'f+z'$. Pairing with $e$ gives
$N(\beta')=D$. The scalar norm equation, without any centering
assumption on this endpoint, gives
$|u'|\ge (\sqrt3/2)(m+Q(z'))/\sqrt D$. By the triangle inequality,
\[
 \sqrt{E(w)}\ge \sqrt{E(n)}\,|u'|-C(\sqrt D+\sqrt{Q(z')})
 \ge \tfrac{\sqrt3}{2}\sqrt{E(n)}\,m/\sqrt D-C'\sqrt D.
\]
For the last step complete the square in $\sqrt{Q(z')}$, using
$E(n)\ge6$. Squaring the positive part, or using
$(x-y)_+^2\ge x^2-2xy$ for $x,y\ge0$, gives
$E(w)\ge(3/4)E(n)m^2/D-O(m)\ge(9/2)m^2/D-O(m)$.
Consequently $E(v)-E(w)=O(m)$ even for thin cusps.

Every strict endpoint decrease spends at least one integer unit.
The bounded-recut comparison and two- and three-slit escape rules take
at most a fixed number of elementary operations between successive
decreases and at termination. Every intermediate tip is obtained
from an endpoint by a bounded word in fixed integral maps, hence has
coordinates $O(\sqrt{E(v)})=O(m)$. This proves the operation and bit-width
claims.
\end{proof}

\begin{corollary}[Output-linear adjacency decompression]
\label{cor:linear-census-adjacency}
Given a selected tuple, its planar adjacency list can be constructed
in deterministic $O(m)$ arithmetic word operations using $O(m)$ words.
The word RAM admits arithmetic with integer division
on $O(\log(2m))$-bit integers. Arbitrary vertex labels suffice.
\end{corollary}
\begin{proof}
Carry out the central translation in the preceding proof by exact
rounding of one rational number. Apply \cref{thm:linear-recut-delay} and the
constant-size ear triangulation of the admissible star. Each fragment
is the positive-area intersection of a chart triangle with one
elementary triangle of the developed mesh. The total
lattice area of its twenty chart triangles is $2m$. All coordinates
are $O(m)$ and the total perimeter, including the nineteen internal
diagonals, is $O(m)$.

Scan each chart by unit-height strips. On each strip the projection
of its convex intersection is one interval. Allocate an array from the
rounded left endpoint minus one to the rounded right endpoint, with
two slots at each integer abscissa for the two lattice orientations.
A unit square allocated this way is contained in the chart plus
$[-3,3]^2$. The interiors of the squares for a single chart are disjoint,
and this enlargement has area $O(1+\text{area}+\text{perimeter})$.
Hence there are $O(m)$ slots altogether, including empty ones.
Compute the endpoints by exact floor and ceiling division of the
unreduced edge intersections. Positive-area overlap of two triangles
uses a fixed number of orientation comparisons. Each occupied slot
stores its fragment index and can be retrieved in constant time.

For each of the thirty seams (nineteen diagonals and eleven paired
slits), transport every source fragment
meeting its interior to the target chart by the seam's lattice
isometry. Test segment intersection with the triangle interior by
separating axes; handle coincidence along an elementary edge by
testing positive interval overlap on that edge. A crossing gives an
edge between equal face fragments, carrying a cyclic slot shift.
Aligned edges give ordinary adjacency links, as do internal chart
neighbors. The target of each transport is found by direct indexing.
Each fragment is examined only a bounded number of times. Traverse
the equal-fragment graph by depth-first search to assign one output
vertex and consistent cyclic offsets to each connected component.
Transport the recorded adjacency links by these offsets and fill
the three slots of each vertex. This yields the same oriented
quotient as seam gluing and costs $O(m)$ word operations, including
the cubic graph and face checks and the writing of all entries.
All orientation products, strip numerators, addresses and counters
are polynomially bounded in $m$.
\end{proof}

\section{Counting completions of a code}
\label{sec:direct-orthoscheme-counts}
Fix the coefficient $\beta$ in the inequalities of \cref{thm:pooled-affine-unranking}. The size equation then has a positive semidefinite
quadratic part in the two caps and a linear part in both scalar
coordinates. Counting this single equation sums all admissible
scalar phases and cap norms together.

\begin{theorem}[A quadratic level with fixed facet directions]
\label{lem:fixed-normal-orthoscheme}
Fix a dimension $d$, a positive semidefinite integral quadratic form $S$
on $\mathbb Z^d$, and a finite set of integral row normals.
Let $P$ be specified by at most one inequality per oriented normal,
with integral right sides of $O(\log(2m))$ bits. Include the coordinate
box $|z_i|\leq C H$ among the inequalities, with $C$ fixed.
For $1\leq H=O(\sqrt m)$, let $L$ be an integral linear form with
coefficients $O(H)$ and let $n$ have $O(\log(2m))$ bits.
The number of points of $P$ with $z\in\mathbb Z^d$ and
$S(z)+L(z)=n$ is computable in
$O(H^2\operatorname{polylog}(2m))$ bit work and bits of storage.
The constants may depend on the form, normals and dimension.
\end{theorem}
\begin{proof}
Replace each right side
$b$ by $b+1/2$. Every admitted lattice point is now an interior
point and every excluded lattice point still violates a row.
Thus an empty or lower-dimensional relaxed polytope contributes
zero. Otherwise enumerate its faces by subsets of the fixed row list
using rational linear algebra in fixed dimension. All faces lie in
the enlarged coordinate box.

For each nonempty face $F$, choose $c_F$ in $\operatorname{aff}F$
minimizing $S$ and then, among its minimizers, minimizing squared
Euclidean length. Writing $V_F$ for the direction space and
$R_F=V_F\cap\ker S$, the defining linear equations say that $c_F$
is $S$-orthogonal to $V_F$ and Euclidean-orthogonal to $R_F$.
They have a unique solution there: on a complement of $R_F$ in
$V_F$ the restriction is positive definite, and the final Euclidean
minimization fixes the radical component. In particular $c_P=0$.
For Euclidean orthogonal projections, Schneider's signed flag formula
gives the polytope indicator away from the simplex walls
\cite[Sec.~3, Proposition~1, Eq.~(5)]{schneider2020volume}.
Here the form can have a radical, and lattice points on shared walls
contribute to the quadratic-level count.
A polytope in its affine
span has a signed decomposition into pyramids from any apex in that
span to its facets. The sign at a facet is positive for an apex on
its inner side, negative on its outer side, and zero on its plane.
Indeed, off the planes of lower-dimensional pieces, follow the ray
from the apex through the test point. Crossings beyond the test point
contribute $+1$ for exit and $-1$ for entry. Their sum is its
indicator. Recursively apply this identity on each facet, with its
own projection as apex. Coning the resulting identity to the previous
apex preserves it away from the cones over the exceptional planes:
the ray meets the facet plane in a unique base point. Consequently
each complete flag contributes the simplex of its projections with
the product of its signs. Discard all zero-volume terms.

Write the flag in decreasing dimension with projections
$c_0=0,c_1,\ldots,c_d$. For consecutive faces with direction
spaces $V_{i-1}\supset V_i$, the difference $c_i-c_{i-1}$
lies in the rational line
$V_{i-1}\cap V_i^{\perp_S}\cap R_i^{\perp_{\rm Eucl}}$,
where $R_i=V_i\cap\ker S$. Indeed on $V_i$ the combined
orthogonality equations have trivial kernel, and impose exactly
$\dim V_i$ independent conditions: use a complement of $R_i$
for the $S$-equations and a basis of $R_i$ for the Euclidean ones.
This proves the line assertion also when the radical dimension changes.
These lines are mutually $S$-orthogonal by nesting. Each face direction is an intersection
of kernels from the fixed normal list. Choose on each line a
primitive integral vector $h_i$, oriented so
$c_i-c_{i-1}=\lambda_i h_i$ with $\lambda_i>0$. There are only
finitely many possible ordered bases $(h_i)$, independently of $H$
and all offsets. Their coordinates and lattice indices are uniformly
bounded. For each fixed direction space the map taking a point in
the affine span to $c_F$ is rational linear, by its defining
equations. Applying it to any face vertex gives $|c_F|=O(H)$,
so $\lambda_i=O(H)$. All lengths are rational
with $O(\log(2m))$-bit descriptions.

In coordinates $z=\sum_i x_i h_i$ the simplex is exactly the chain
\[
       1\geq x_1/\lambda_1\geq x_2/\lambda_2\geq
                 \cdots\geq x_d/\lambda_d\geq0.
\]
To assign all boundaries, evaluate only the polyhedral indicators
at $z+(\delta,\delta^2,\ldots,\delta^d)$ for a formal positive
infinitesimal $\delta$. Every nonzero affine wall has a definite
side there. The relaxed input walls keep their values at $z$.
Each chain comparison becomes either weak or strict on the
unperturbed coordinates, according to the first nonzero perturbation
coefficient. Thus the signed identity holds at every lattice point,
including those on flags shared by several simplices. Throughout
the energy $S+L$ is evaluated at $z$ itself.

Split the lattice into cosets of $\sum_i\mathbb Z h_i$. In one coset
the $x_i$ range independently over $\rho_i+\mathbb Z$ for a fixed
rational shift vector $\rho$. Clear its fixed denominators by a
positive integer $D$ so that each
$e_i(x_i)=D(S(h_i)x_i^2+L(h_i)x_i)$ is integral. The sum of the
unary energies is $D(S(z)+L(z))$ and may be negative.

Form one event for each possible value $0\leq x_i\leq\lambda_i$
on each axis. Enforce the two outer endpoint conditions on their
respective axes. Sort the $O(H)$ events by decreasing $x_i/\lambda_i$.
For equal ratios, order label $i$ before $i+1$ if that adjacent
comparison admits equality, and after it otherwise. A single label
order satisfying these prescriptions exists because they orient a
path. Thus a subsequence with labels $1,\dots,d$ is exactly a
permitted chain. Different coordinate tuples give different
subsequences even at coincident ratios.

Bisect the event list into a balanced binary tree. For each node and
each consecutive label block $[i,j]$ store the Laurent polynomial
counting its allowed subsequences by total energy. The empty block
has polynomial one. At a leaf a matching singleton contributes its
one monomial. At an internal node sum over $k=i-1,\dots,j$ the
product of the left polynomial for $[i,k]$ and the right one for
$[k+1,j]$. Each subsequence splits uniquely between the children.
The answer for the simplex and coset is the coefficient at $Dn$
for $[1,d]$ at the root.

For a node $J$ the support width of a nonzero block polynomial is
at most the sum, over its labels, of the energy spans among their
events in $J$. Store with an exponent offset, so a high-energy
monomial needs only one coefficient. On each tree level the
successive event intervals for any one label are disjoint in its
sorted coordinate list. Their energy spans sum to at most the
total variation of the unary quadratic on $[0,\lambda_i]$,
which is $O(H^2)$ since its linear coefficient is $O(H)$.
Adding one slot per node and summing over the fixed number of
blocks leaves $O(H^2)$ slots per level. Each product and its
alignment in the parent can be charged to that parent's span.
Coefficients count at most $(O(H)+1)^d$ tuples. Packing coefficients
in sufficiently wide $O(\log(2m))$-bit digits and using fast
integer multiplication \cite{harveyHoeven2021} gives quasilinear bit work in the number
of slots up to logarithmic factors in $2m$, with exact carries.
Sum over tree levels, flags and cosets to obtain the claimed bounds.
\end{proof}

\begin{corollary}[The complete distribution of quadratic levels]
Under the hypotheses of \cref{lem:fixed-normal-orthoscheme}, the table
\[
 A(t)=\#\{z\in P\cap\ZZ^d:S(z)+L(z)=t\}
\]
for all attainable integer levels $t$ can be written in
$O(H^2\operatorname{polylog}(2m))$ bit work and storage. An explicit
table can require $\Omega(H^2)$ entries even in dimension two.
\end{corollary}
\begin{proof}
The root polynomial for each flag and coset already contains the
coefficients for every level. For each integer $t$, take its coefficient
at the exponent $Dt$ for that coset and add with the flag sign. There
are a fixed number of such sums per level. Since $|z_i|\le CH$ and
the coefficients of $L$ are $O(H)$, all levels lie in an integer
interval of length $O(H^2)$. For sharpness let $H$ be a positive
integer, $P=[0,H-1]^2$, $S=0$ and $L(x,y)=x+Hy$. Each of the
$H^2$ integers from $0$ to $H^2-1$ occurs once.
\end{proof}

\begin{corollary}[Random access to a quadratic level]
Under the hypotheses of \cref{lem:fixed-normal-orthoscheme}, one can
rank a lattice point or unrank an integer in lexicographic order in
$O(H^2\operatorname{polylog}(2m))$ bit work and storage.
\end{corollary}
\begin{proof}
Adjoin the $2d$ signed coordinate normals to the fixed list. Fixing a
prefix uses pairs of opposite inequalities. An upper bound on the next
coordinate gives its cumulative count. Choose the first threshold whose
count exceeds the zero-based residual rank and subtract the count at the
preceding threshold. For ranking, sum those preceding counts along the
given point. At most $O(d\log(2H))$ counts suffice, reusing the workspace.
\end{proof}

For example, consider $0\leq y<x\leq5$ with energy $x^2+y^2$.
This triangle already has orthogonal coordinate axes. At an equal
ratio the event for $y$ precedes the event for $x$, so the pair
$(3,3)$ contributes nothing to the coefficient at $18$. At energy
$25$ the subsequences for $(4,3)$ and $(5,0)$ both survive, including
the point on the two closed outer edges. In the general signed sum
the infinitesimal rule makes exactly these boundary decisions for
each simplex and for every shared face.

\begin{theorem}[Direct all-rank selection by counting family completions]
\label{thm:orthoscheme-family-count}
For every $N=2m$ and every coordinate prefix and next-coordinate
threshold in the order $(H,a,z_1,\ldots,z_{16},s)$, the number of
canonical completions is computable in deterministic
$O(m^2\operatorname{polylog}(2m))$ bit operations and
$O(m\operatorname{polylog}(2m))$ bits of working storage.
These bounds include setup, also for the empty prefix. For every
integer rank from zero through the isomer count minus one, they
suffice to construct the selected graph from the size and rank.
The algorithm visits height coefficients. Its tables count partial
energies without listing scalar phases, individual caps or cap pairs.
\end{theorem}
\begin{proof}
Fix $\beta=a+b\zeta$, $H=a+b$, and set $p_0=\lfloor m/H\rfloor$.
For $b>0$ set $q_0=0$ and for $b=0$ set $q_0=-\lfloor p_0/2\rfloor$.
Put $P'=p-p_0$ and $T'=q-q_0$. On every accepted tuple
$|z_i|\leq H$ and $Q(z)=O(H^2)$. For $b>0$ the strip gives
$|q|\leq H$, and $Hp=m+Q(z)-bq$ then gives $|P'|=O(H)$.
For $b=0$ the same equation gives $|P'|=O(H)$ and
$|p+2q|\leq H$ gives $|T'|=O(H)$. Thus in both cases we may
impose a fixed box of radius $CH$ in $(z,P',T')$. The size test is
\[
 Q(z)-HP'-bT'=Hp_0+bq_0-m .
\]
The quadratic part is fixed. Its radical consists of the two scalar
directions, and its restriction to the cap coordinates is positive definite.
In the predicate of \cref{thm:pooled-affine-unranking}
every ownership wall is a fixed linear form in $(z,P',T')$ with an integral
threshold after clearing fixed denominators. The coordinate box and prefix constraints add only fixed normals,
since cap cuts fix or bound individual integral coordinates and the final scalar cut uses
$s=q$ or $s=p+2q$ on the two respective branches. Before the coefficient
is fixed, a prefix or threshold in $(H,a)$ simply selects which coefficient
counts enter the sum.
Remove constant tests directly. Form disjoint sign cells, converting
strict signs to integer thresholds and equality to two inequalities
before relaxing any cells. The number of rows per cell is bounded
independently of $m$. For repeated oriented normals retain their
tightest bound. Each integer owner retains exactly one
cell, including on tie and mirror boundaries. Apply
\cref{lem:fixed-normal-orthoscheme} to the displayed size test on each cell.

This costs $O(H^2\operatorname{polylog}(2m))$ per $\beta$.
There are $O(H)$ preferred coefficients at height $H$.
Summing $O(H^3)$ up to $\sqrt{54m}$ gives the time bound. Reusing
the tables between coefficients gives the space bound. All nineteen code
coordinates have polynomial ranges. Their cumulative counts select
each in turn by binary search, taking $O(\log(2m))$ calls in total.
Recover the scalar by the size equation and strip coordinate, center
the bounded cap tuple by the scalar centering described in \cref{thm:linear-recut-delay} and apply \cref{thm:linear-recut-delay,cor:linear-census-adjacency}.
The graph uses $O(m)$ additional words and $O(m)$ word operations.
\end{proof}

\begin{corollary}[A cumulative index by coefficient]
\label{cor:orthoscheme-height-directory}
In $O(m^2\operatorname{polylog}(2m))$ bit work one may build an index
of $O(m\log(2m))$ bits from which each rank returns its graph in
$O(m\operatorname{polylog}(2m))$ bit work. Construction and lookup
each use $O(m\operatorname{polylog}(2m))$ bits of working storage.
\end{corollary}
\begin{proof}
Compute the mass of each preferred $\beta=(a,b)$ once by the theorem's
fixed-coefficient algorithm. Store cumulative masses in lexicographic
$(H,a)$ order, the first two entries of the code. There are
$O(m)$ entries, each of logarithmic bit length. For a zero-based rank
$r$, the first mass total strictly exceeding $r$ determines $\beta$;
subtract the preceding total. Repeated totals from empty families are
therefore harmless. With $\beta$ now fixed, each remaining bisection
query costs $O(H^2\operatorname{polylog}(2m))$. Since $H^2<54m$,
the $O(\log(2m))$ queries and graph decoding fit within the stated
bound. Only the cumulative counts and coefficient keys persist.
\end{proof}

\begin{proof}[Proof of \cref{thm:main}]
The period correspondence and \cref{thm:pooled-affine-unranking}
give a bijection from the codes to ordinary isomers. The last entry
of the cumulative index in \cref{cor:orthoscheme-height-directory}
is their number (zero for an empty index). For a rank in the stated
range, bisection gives one complete code. Its norm equation recovers
the scalar, and \cref{thm:star-complete-recut-normalization} recovers
the original period orbit in the output graph. The indexed bit bounds
and the word bound now follow from
\cref{cor:orthoscheme-height-directory,cor:linear-census-adjacency}.
\end{proof}

\subsection{Ranking arithmetic codes and sampling ordinary isomers}
\begin{corollary}[Rank intervals of arithmetic codes]\label{cor:rank-interval}
With the index of \cref{cor:orthoscheme-height-directory}, a prefix
of the nineteen code coordinates determines a half-open interval
$[L,U)$ of zero-based ranks in $O(m\operatorname{polylog}(2m))$
bit operations. Its length is the number of completions. For a
complete admissible code, $L$ is its rank.
\end{corollary}
\begin{proof}
Lexicographic order makes the completions consecutive. For a prefix
of length at most two, two binary searches on the sorted coefficient
keys give the boundaries, and differences of cumulative masses give
their ranks. The empty prefix gives $[0,\iso(2m))$. For a longer
prefix, first look up its key $(H,a)$, which fixes $b=H-a$. An absent pair contributes
an empty interval at its insertion rank. For a present pair, start
$L$ at its preceding cumulative mass. At each specified remaining
coordinate, add the number of completions agreeing with all preceding
coordinates and strictly smaller at this coordinate. These disjoint
sets partition the words of this pair below the prefix. Finally add
the number agreeing with the whole prefix to obtain $U$ from $L$.
Each count uses the fixed-coefficient calculation in
\cref{thm:orthoscheme-family-count}, at cost
$O(H^2\operatorname{polylog}(2m))$. There are only seventeen
remaining coordinates and $H^2<54m$. The singleton assertion follows
from uniqueness of the representative. The ranking input here is the
arithmetic code of the isomer.
\end{proof}

\begin{corollary}[Exact uniform sampling from the vertex number]
\label{cor:uniform-sampling}
Given $2m$, the deterministic preprocessing of \cref{thm:main}
detects an empty family. Otherwise each subsequent independent draw
returns an ordinary isomorphism class with probability $1/\iso(2m)$
in expected $O(m\operatorname{polylog}(2m))$ bit operations.
The same bound holds for the uniform distribution conditional on any
nonempty arithmetic prefix, including the cost of finding its interval.
\end{corollary}
\begin{proof}
Write $M=U-L$ for the interval length from \cref{cor:rank-interval}.
If $M=0$ report the empty family. If $M=1$ decode $L$. Otherwise draw
$k=\lceil\log_2 M\rceil$ independent fair bits as an integer in
$[0,2^k)$. Repeat with fresh bits until the integer $t$ is below $M$,
and unrank $L+t$. Each accepted value has probability $1/M$, since
every trial assigns the same probability $2^{-k}$ to it. The expected
number of trials is $2^k/M<2$ and $k=O(\log(2m))$ by the polynomial
coordinate ranges. Only one graph is decoded per draw. Because each
ordinary isomer has one rank, automorphism orders do not enter the
sampling probabilities. Fresh bits give independent draws from the
same precomputed index. For example, at $2m=60$ the full interval has
length $1812$ \cite[Sec.~3.3]{bille2025open}. An eleven-bit trial is
accepted with probability $1812/2048$, and the conditional probability
of each isomer is $1/1812$.
\end{proof}

\begin{figure}[tbp]
\centering
\includegraphics[width=.93\linewidth]{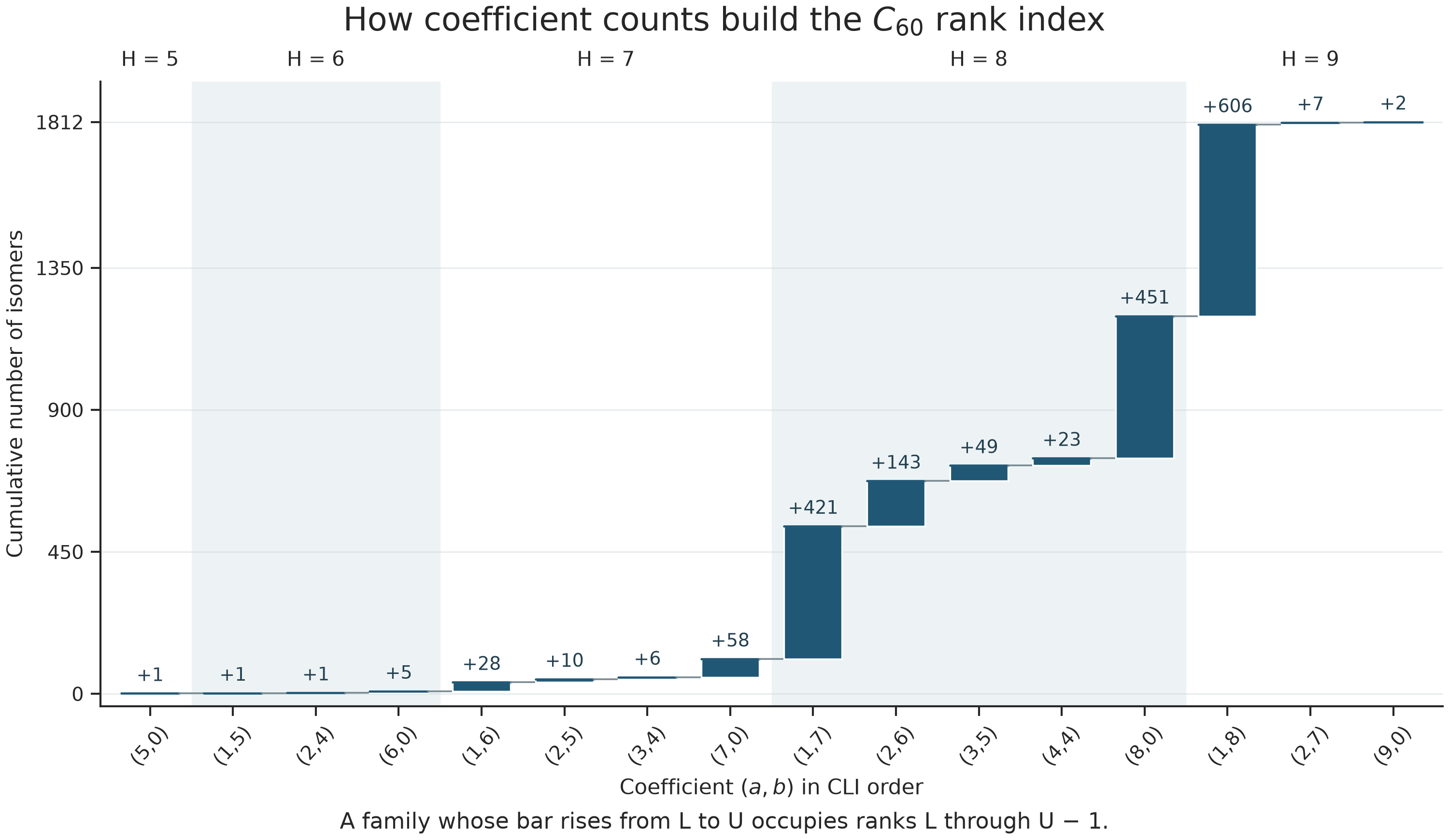}
\caption{The coefficient index for $\mathrm C_{60}$. The $16$ nonempty
families appear in $(H,a)$ order. A bar starts at the sum of all
preceding family counts and rises by the count printed above it.
For example, $(1,8)$ contributes $606$ isomers, with ranks $1197$
through $1802$.}
\label{fig:c60-coefficient-index}
\end{figure}

\paragraph{Unranking zero-based rank nine of $\mathrm C_{60}$.}
The first five bars in \cref{fig:c60-coefficient-index} locate the coefficient.
Take $m=30$ and $r=9$. The nonempty coefficient families up to
$(a,b)=(1,6)$ and their half-open rank intervals are
\[
\begin{array}{c|ccccc}
(a,b)&(5,0)&(1,5)&(2,4)&(6,0)&(1,6)\\
\hline
[L,U)&[0,1)&[1,2)&[2,3)&[3,8)&[8,36).
\end{array}
\]
Subtracting eight gives local rank one in the $(1,6)$ family.
In increasing cap-pair order and then increasing $s=q$, the
selected caps are
\begin{align*}
 z_A&=(-7,-7,\;-7,-6,\;-7,-7,\;-7,-3),\\
 z_B&=(-7,-7,\;-7,-6,\;-7,-7,\;-7,-2).
\end{align*}
Consecutive pairs give the two integer coordinates of an Eisenstein
entry. At the preceding rank the cap pair is $(z_A,z_A)$ with $q=-1$.
At rank nine the second cap increases, so its scalar can be smaller. The scalar is $q=-3$ and the cap energy is $274$, so
$7p+6(-3)=30+274$ gives $p=46$.
Decoding this period gives a cubic graph on 60 vertices. Tracing the
faces in its cyclic adjacency lists gives twelve pentagons and twenty
hexagons. Thirteen edges are shared by two pentagons.
The commands
\begin{verbatim}
buckycount unrank 60 9 -o c60-r9.pc
buckycount unrank 60 1 --family 1,6 -o c60-family.pc
\end{verbatim}
write identical planar-code files for this isomer.

\section{Conclusion}
Within a fixed set of facet directions the counting algorithm returns
the whole distribution of a quadratic statistic at once. Its exponent
in the box radius is sharp for an explicit table. For random access
one instead uses a single level and repeats the count after coordinate cuts.
For fullerenes, two geometric results make this useful: each isomer
has one arithmetic representative selected by fixed-direction walls,
and its period can be decoded into an adjacency list in linear word time.
For other orbit problems, the first question is whether their
representatives admit such a description with uniform coefficient bounds.

\section*{Declarations}
\subsection*{Use of artificial intelligence}
During preparation of this work, the author used Fable 5.2 (Claude Code) and
Astra 6.2 (OpenAI Codex) and Deep Review (Google) to assist with exploration,
manuscript revision and verification code.  After using these tools, the
author reviewed and edited the content as needed and takes full
responsibility for the content of the publication.
\subsection*{Competing interests}
The author declares no conflict of interest.
\subsection*{Data and code availability}
The archive DOI for the fullerene counting paper series and its code and data packages is
\url{https://doi.org/10.5281/zenodo.24289727}.
The ancillary package includes the source of \texttt{buckycount}, which counts and unranks fullerenes, and the data for the star recuts.

\bibliographystyle{plain}
\bibliography{references}
\end{document}